\documentclass[11pt]{amsart}
\usepackage{amsfonts}
\usepackage{mathrsfs}
\usepackage{amssymb}
\usepackage{hyperref}
\usepackage{amsmath,amssymb, amsthm}
\usepackage{tikz}
\usetikzlibrary{arrows,shapes}

\newcommand{\SL}{\mathrm{SL}(2,\mathbb{R})}

\newcommand{\rot}{\mathrm{rot}}

\let\newpf\proof \let\proof\relax 
\newenvironment{pf}{\newpf[\proofname]}{\qed\endtrivlist}

\newcommand{\ba}{\overline{A}}

\def\be{\begin{equation}}
\def\ee{\end{equation}}

\def\ba{{\begin{align}}}
\def\ea{{\end{align}}}

\def\bm{\begin{matrix}}
\def\em{\end{matrix}}

\def\u{{\mathbb U}}

\def\x{{\bf x}}

\def\SL{{\mathrm{SL}}}

\def\0{{\mathbf 0}}

\newtheorem{Theorem}{Theorem}[section]
\newtheorem{Lemma}{Lemma}[section]
\newtheorem{Proposition}{Proposition}[section]

\newtheorem{Remark}{Remark}[section]

\newtheorem{Definition}{Definition}[section]

\numberwithin{equation}{section}

\theoremstyle{definition}

\def\ssm{\smallsetminus}
\def\tr{{\text{tr}}}
\renewcommand{\setminus}{\ssm}

\renewcommand{\mod}{\operatorname{mod}}

\newcommand{\id}{\operatorname{id}}

\newcommand{\C}{{\mathbb C}}

\newcommand{\N}{{\mathbb N}}
\newcommand{\Q}{{\mathbb Q}}
\newcommand{\R}{{\mathbb R}}
\newcommand{\T}{{\mathbb T}}

\newcommand{\Z}{{\mathbb Z}}

\def\B0{{\bold{0}}}

\catcode`\@=12

\def\Empty{}
\newcommand\oplabel[1]{
  \def\OpArg{#1} \ifx \OpArg\Empty {} \else
    \label{#1}
  \fi}

\newcommand{\comm}[1]{}
\newcommand{\comment}[1]{}

\begin{document}

\title[Phase transition]{Sharp Phase transitions  for the almost Mathieu operator}

\author{Artur Avila}
\address{CNRS UMR 7586, Institut de Math\'ematiques de Jussieu -
Paris Rive Gauche, B\^atiment Sophie Germain, Case 7012, 75205 Paris
Cedex 13, France \& IMPA, Estrada Dona Castorina 110, 22460-320, Rio
de Janeiro, Brazil}
 \email{artur@math.jussieu.fr}

\author {Jiangong You}
\address{
Department of Mathematics, Nanjing University, Nanjing 210093, China
} \email{jyou@nju.edu.cn}

\author{Qi Zhou}
\address{Laboratoire de Probabilit\'{e}s et Mod\`{e}les
al\'{e}atoires, Universit\'{e} Pierre et Marie Curie, Boite courrier
188 75252, Paris Cedex 05, France
}
\curraddr{Department of Mathematics, Nanjing University, Nanjing 210093, China}

 \email{qizhou628@gmail.com}

\date{\today}

\begin{abstract}
It is known that the spectral type of the almost Mathieu operator depends in a fundamental way on both the strength of the coupling constant and the arithmetic properties of the frequency.  We study the competition between those factors and locate the point where the phase transition from singular continuous spectrum to pure point spectrum takes place, which
 solves Jitomirskaya's conjecture in
\cite{Ji95,J07}. Together with \cite{Aab}, we give the  sharp
description of phase transitions for the almost Mathieu operator.
\end{abstract}

\setcounter{tocdepth}{1}

\maketitle

\section{Main results}

This paper concerns the spectral measure of the Almost Mathieu operator:
\begin{equation*}
(H_{\lambda,\alpha,\theta} u)_n= u_{n+1}+u_{n-1} +2\lambda \cos 2
\pi (n\alpha + \theta) u_n,
\end{equation*}
where $\theta\in \mathbb{R}$ is  the phase, $\alpha\in \R\backslash
\Q$ is the frequency and $\lambda\in \R$ is  the coupling constant,
which has been extensively  studied because of its strong backgrounds in
physics and also because it provides interesting examples in spectral theory
\cite{L1}. We will find  the exact transition point from singular
continuous spectrum to purely point spectrum of Almost Mathieu operator, thus solve  Jitomirskaya's conjecture in 1995 \cite{Ji95}(see also Problem 8 in \cite{J07}).
More precisely, let $\frac{p_n}{q_n}$ be the $n-$th convergent of $\alpha,$ and
define
\begin{equation}\label{defbeta}\beta(\alpha):=\limsup_{n\rightarrow \infty}\frac{\ln
q_{n+1}}{q_n},\end{equation} our main results are the following:

\begin{Theorem}\label{main}
Let $\alpha\in \R\backslash \Q$ with $0<\beta(\alpha)<\infty$, then
we have the following:
\begin{enumerate}
\item If $|\lambda|<1,$ then $H_{\lambda,\alpha,\theta}$ has
purely absolutely continuous spectrum for all $\theta$.
\item  If $1\leq |\lambda|<e^\beta,$ then $H_{\lambda,\alpha,\theta}$ has
purely singular continuous spectrum for all $\theta$.
\item  If $|\lambda|>e^\beta,$ then $H_{\lambda,\alpha,\theta}$ has
purely point spectrum with exponentially decaying eigenfunctions
for a.e. $\theta$.
\end{enumerate}

\end{Theorem}

\begin{Remark}
 Part (1) is proved by Avila \cite{Aab}, we state here just for
completeness.\end{Remark}
\begin{Remark} The cases $\beta=0, \infty$ have been solved in previous works \cite{Aab,AJ05}. Together with Theorem \ref{main},
 one sees the sharp phase transition scenario of three types of the spectral measure. Moreover, the type of the spectral measure
 is clear for all $(\lambda, \beta)$ except the line $\lambda=e^\beta$. See Figure 1 and Figure 2 below.
\end{Remark}

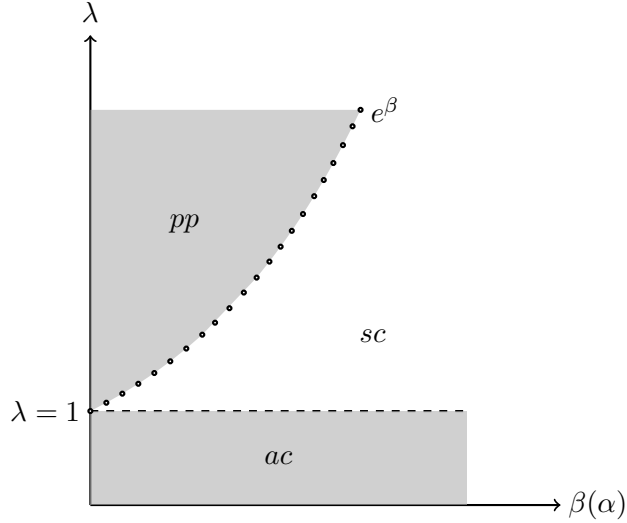
\begin{figure}[th]
  \centering
  \begin{tikzpicture}[scale = 2.5,x={(1cm,0cm)},y={(0cm,0.5cm)}]
    \tikzset{mypoints/.style={fill=white,draw=black,thick}}

    \coordinate (origin) at (0,0);
    \coordinate (beta) at (2.5,0);
    \coordinate (lambda) at (0,5);

    \draw [draw,->, thick] (origin) -- (beta);
    \draw [draw,->, thick] (origin) -- (lambda);

    \coordinate (ac_1) at (0,1);
    \coordinate (ac_2) at (2,1);
    \coordinate (ac_3) at (2,0);
    \draw [thick,dashed] (ac_1) -- (ac_2);

    \fill[gray!60,opacity=0.6] (origin)--(ac_1)--(ac_2)--(ac_3)--cycle;

    \coordinate[label = left: {$\lambda = 1$}] (A) at (0*0.17,{(e)^(0*0.17)});
    \coordinate (B) at (0.5*0.17,{(e)^(0.5*0.17)});
    \coordinate (C) at (1*0.17,{(e)^(1*0.17)});
    \coordinate (D) at (1.5*0.17,{(e)^(1.5*0.17)});
    \coordinate (E) at (2*0.17,{(e)^(2*0.17)});
    \coordinate (F) at (2.5*0.17,{(e)^(2.5*0.17)});
    \coordinate (G) at (3*0.17,{(e)^(3*0.17)});
    \coordinate (H) at (3.5*0.17,{(e)^(3.5*0.17)});
    \coordinate (I) at (3.9*0.17,{(e)^(3.9*0.17)});
    \coordinate (J) at (4.35*0.17,{(e)^(4.35*0.17)});
    \coordinate (K) at (4.8*0.17,{(e)^(4.8*0.17)});
    \coordinate (L) at (5.2*0.17,{(e)^(5.2*0.17)});
    \coordinate (M) at (5.6*0.17,{(e)^(5.6*0.17)});
    \coordinate (N) at (5.95*0.17,{(e)^(5.95*0.17)});
    \coordinate (O) at (6.3*0.17,{(e)^(6.3*0.17)});
    \coordinate (P) at (6.65*0.17,{(e)^(6.65*0.17)});
    \coordinate (Q) at (7*0.17,{(e)^(7*0.17)});
    \coordinate (R) at (7.3*0.17,{(e)^(7.3*0.17)});
    \coordinate (S) at (7.6*0.17,{(e)^(7.6*0.17)});
    \coordinate (T) at (7.9*0.17,{(e)^(7.9*0.17)});
    \coordinate (U) at (8.2*0.17,{(e)^(8.2*0.17)});
    \coordinate[label = right: {$e^{\beta}$}] (V) at (8.45*0.17,{(e)^(8.45*0.17)});
    \coordinate (OO) at (0,{(e)^(8.45*0.17)});

    \fill[gray!60,opacity=0.6] (A)--(B)--(C)--(D)--(E)--(F)--(G)--(H)--(L)--(M)--(N)--(O)--(P)--(Q)--(R)--(S)--(T)--(U)--(V)--(OO)--cycle;

    \foreach \p in {A,B,C,D,E,F,G,H,I,J,K,L,M,N,O,P,Q,R,S,T,U,V}
        \fill[mypoints] (\p) circle (0.3pt);

    \draw(2.7,0) node{$\beta(\alpha)$};
    \draw(0,5.25) node{$\lambda$};

    \draw(1,0.5) node{$ac$};
    \draw(1.5,1.8) node{$sc$};
    \draw(0.5,3) node{$pp$};
  \end{tikzpicture}
  \caption{Phase transition diagram}
  \label{fig:1}
\end{figure}

\begin{figure}[th]
  \centering
  \begin{tikzpicture}[scale = 2.5,x={(1cm,0cm)},y={(0cm,0.5cm)}]
    \coordinate (origin) at (0,0);
    \coordinate (A) at (1,0);
    \coordinate (B) at (2.5,0);
    \coordinate[label = right: {$\lambda$}] (C) at (4.2,0);

    \draw(0.5,0.2) node{$ac$};
    \draw(1.75,0.2) node{$sc$};
    \draw(3.25,0.2) node{$pp$};

    \draw [draw,->, thick] (origin) -- (C);

    \foreach \x/\xtext in {0,1}
      \draw[xshift=\x cm] (0pt,1pt) -- (0pt,-1pt) node[below,fill=white]
            {$\xtext$};
    \draw[xshift=2.5 cm] (0pt,1pt) -- (0pt,-1pt) node[below,fill=white]
            {$e^{\beta}$};

  \end{tikzpicture}
  \caption{Phase transition for fixed $\alpha$.}
  \label{fig:2}
\end{figure}
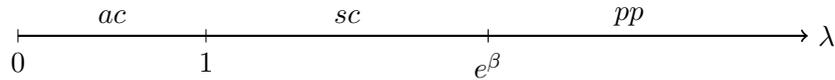

\begin{Remark} 
Theorem \ref{main}(3), also called Anderson localization (AL),  is optimal in the sense
that the result can not be true for $G_{\delta}$ dense $\theta$  \cite{JS}. The
arithmetic property of $\theta$ will influence the spectral measure. \end{Remark}

 Now we
briefly recall the history of this problem. By symmetry, we just need to consider the case $\lambda>0$. In 1980, Aubry-Andr\'e
\cite{AA80} conjectured that the spectral measure of
$H_{\lambda,\alpha,\theta}$ depends on $\lambda$ in the following
way:
\begin{enumerate}
\item If $\lambda<1,$ then $H_{\lambda,\alpha,\theta}$ has
purely absolutely continuous spectrum for all $\alpha\in
\R\backslash \Q$, and  all $\theta\in\R$.
\item   If $\lambda >1,$ then $H_{\lambda,\alpha,\theta}$ has
pure point spectrum for all $\alpha\in \R\backslash \Q$, and  all
$\theta\in \R$.
\end{enumerate}

However,  Aubry and Andr\'e overlooked the role of the arithmetic
property of $\alpha$.   Avron-Simon \cite{AS} soon found  that by
Gordon's lemma \cite{G},  $H_{\lambda,\alpha,\theta}$ has no
eigenvalues for  any $\lambda\in\R$, $\theta\in\R$ if $\beta(\alpha)=\infty$. Since then,
people pondered how the arithmetic property of $\alpha$ influences
the spectral type and under which condition Aubry-Andr\'e's
conjecture \cite{AA80} is  true.

When $\alpha$ is Diophantine (i.e. there exist $\gamma,\tau>0$  such that 
$\|k\alpha\|_{\T} \geq
\frac{\gamma^{-1}}{|k|^{\tau}},$ for all $0
\neq k \in \Z$), and $\lambda$ is large enough, the operator has
pure point spectrum \cite{E97,FSW,Sin}, and when $\lambda$ is
small enough, the operator has absolutely continuous spectrum
\cite{CD,DS75,E92}. The common feature of the above results is that
they both rely on KAM-type arguments, thus the largeness or
smallness of $\lambda $ depend on Diophantine constant
$\gamma,\tau$,  we therefore call such results  perturbative results.
Non-perturbative approach to localization  problem was
developed by Jitomirskaya, based on partial advance \cite{J94,J95}, she 
finally  proved that if  $\alpha$ is Diophantine, 
$H_{\lambda,\alpha,\theta}$ has AL for all  $\lambda>1$ and a.e. $\theta\in\R$. It
follows from  the strong version of Aubry duality \cite{GJLS}, 
$H_{\lambda^{-1},\alpha,\theta}$ has purely absolutely continuous
spectrum for a.e. $\theta\in\R$. Therefore, Jitomirskaya \cite{J99}
proved Aubry-Andr\'e's conjecture in the measure setting, i.e. the
conjecture holds for almost every $\alpha\in \R\backslash \Q,$
$\theta\in \mathbb{R}$.

Before Jitomirskaya's result, Last \cite{L93}, Gesztesy-Simon \cite{GS},
 Last-Simon \cite{LS}
have already  showed that $H_{\lambda,\alpha,\theta}$ has absolutely
continuous components for every $\lambda<1$, $\alpha\in \R\backslash
\Q,$ $\theta\in \mathbb{R}$, so the conjecture in subcritical regime
still has some hope to be true, which was  also conjectured by Simon
\cite{Si00}. Recently, Avila-Jitomirskaya \cite{AJ08} showed that if
$\alpha$ is Diophantine, then $H_{\lambda,\alpha,\theta}$ is purely
absolutely continuous for every $\theta\in \mathbb{R}$. For
$\beta>0$,  Avila-Damanik \cite{AD} proved that the conjecture (1)
 for almost every $\theta$. The complete answer of  
Aubry-Andr\'e's conjecture (1) was provided by  Avila \cite{Aab}. One thus sees 
that $\lambda=1$ is the phase transition point from absolutely continuous spectrum to singular spectrum.

The remained issue is Aubry-Andr\'e's conjecture (2) when
$\alpha$ is  Liouvillean.  People already knew that the spectral measure is pure point for Diophantine $\alpha$ and almost every phases, while it is purely singular continuous for $\beta(\alpha)=\infty$ and all phase. So there must be phase transition when $\beta(\alpha)$ goes from zero to infinity.  In 1995, Jitomirskaya \cite{Ji95} 
modified the second part of the Aubry-Andr\'e's conjecture and conjectured the following\begin{enumerate}
\item  If $1<\lambda<e^\beta$, the spectrum is purely singular continuous for all $\theta$.
\item  If $\lambda>e^\beta$, the spectrum is pure
point with exponential decaying eigenfunctions for a.e. $\theta$.
\end{enumerate}
 Thus $\lambda=e^\beta$ is conjectured to be
 the exact phase transition point from  continuous spectrum to pure point spectrum. There are some partial results on Jitomirskaya's conjecture. By Gordon's lemma \cite{G} and the
exact formula of Lyapunov exponent \cite{BJ}, one can prove that
$H_{\lambda, \alpha,\theta}$ has purely singular continuous spectrum
for any $\theta\in\R$ if $1<\lambda<e^{\frac{\beta}{2}}$, see also
Remark \ref{gordon} for more discussions. For the pure point part,
Avila-Jitomirskaya \cite{AJ05} showed that if
$\lambda>e^{\frac{16\beta}{9}}$, then $H_{\lambda,\alpha,\theta}$
has AL for a.e. $\theta\in\R$.  You-Zhou \cite{YZ} proved that if $\lambda>Ce^{\beta}$ with $C$ large enough \footnote{If one check carefully the proof, it already gives $C=1$.}, then the eigenvalues of  $H_{\lambda,\alpha,\theta}$ with exponentially decaying eigenfunctions are dense in the spectrum.  Readers can find more discussions on these two results  in section \ref{ander}.
The  main contribution of this paper is to give a full  proof of Jitomirskaya's conjecture. 

We remark that the spectral type at the transition points $\lambda=1$ and $\lambda=
e^\beta$ have not been completely understood so far. Partial results include the
following: in case $\lambda=1$, since the Lebesgue measure of the
spectrum is zero for every $\alpha\in \R\backslash \Q$
\cite{AK06,L94},  by Aubry duality \cite{GJLS}, we know $H_{\lambda,
\alpha,\theta}$ is purely singular continuous for a.e.
$\theta\in\R$. In fact, Avila \cite{App} has proved more: if
$\theta$ is not rational w.r.t $\alpha$, then $H_{\lambda,
\alpha,\theta}$ is purely singular continuous. We remark that, by
Gordon's lemma \cite{G}, if $\beta > 0$, then $H_{\lambda,
\alpha,\theta}$ is purely  singular continuous for $\lambda=1$ and
every $\theta\in\R$, we  include this  in Theorem \ref{main}(2).
Excluding or proving the existence of point spectrum in case that
$\alpha $ is Diophantine is one of the major interesting problems for
the critical almost Mathieu operator. For the second transition
point $\lambda=e^\beta$, one knows almost nothing but purely
singular continuous spectrum  for a $G_{\delta}$ set of $\theta$
\cite{JS}. The spectral type possibly depends on the finer
properties of approximation of $\alpha$, as conjectured by
Jitomirskaya in \cite{J07}.

\section{preliminaries}

 For a bounded
analytic (possibly matrix valued) function $F$ defined on $ \{ \theta |  | \Im \theta |< h \}$, let
$
\|F\|_h=  \sup_{ | \Im \theta |< h } \| F(\theta)\| $ and denote by $C^\omega_{h}(\T,*)$ the
set of all these $*$-valued functions ($*$ will usually denote $\R$,
$SL(2,\R)$). 

\subsection{Continued Fraction Expansion}\label{sec:2.1}
Let $\alpha \in (0,1)$ be irrational. Define $ a_0=0,
\alpha_{0}=\alpha,$ and inductively for $k\geq 1$,
$$a_k=[\alpha_{k-1}^{-1}],\qquad \alpha_k=\alpha_{k-1}^{-1}-a_k=G(\alpha_{k-1})=\{{1\over \alpha_{k-1}}\},$$
Let $p_0=0,  p_1=1,  q_0=1,  q_1=a_1,$ then we define inductively
$p_k=a_kp_{k-1}+p_{k-2}$, $q_k=a_kq_{k-1}+q_{k-2}.$
The sequence $(q_n)$  is the  denominators of best rational
approximations of $\alpha$ since we have \begin{equation} \forall 1
\leq k < q_n,\quad \|k\alpha\|_{\T} \geq \|q_{n-1}\alpha\|_{\T},
\end{equation}
and
\begin{equation}
\|q_n \alpha \|_{\T} \leq {1 \over q_{n+1}}.
\end{equation}
Note that $(\ref{defbeta})$ is equivalent to 
\be\label{equibeta}
\limsup_{k\rightarrow \infty} \frac{1}{|k|} \ln \frac{1}{ \|k\alpha\|_{\T}}=\beta.
\ee

\subsection{Cocycles} A cocycle $(\alpha, A)\in \R\backslash
\Q\times C^\omega(\T, SL(2,\R))$ is a linear skew product:
\begin{eqnarray*}\label{cocycle}
(\alpha,A):&\T^{1} \times \R^2 \to \T^{1} \times \R^2\\
\nonumber &(\theta,v) \mapsto (\theta+\alpha,A(\theta) \cdot v),
\end{eqnarray*}
for $n \geq 1$, the products are defined as
$$A_n(\theta)=A(\theta+(n-1)\alpha) \cdots
A(\theta),$$ and $A_{-n}(\theta)=A_n(\theta-n\alpha)^{-1}.$ For this kind of cocycles, the Lyapunov exponent  
$$ L(\alpha, A)=\lim_{n\rightarrow \infty} \frac {1} {n}
\int \ln \|A_n(\theta)\| d\theta,
$$
is well defined.

Assume now  $A\in C^0(\T, SL(2,\R))$ is homotopic to the
identity. Then there exists $\psi:\T \times \T \to \R$ and $u:\T
\times
\T \to \R^+$ such that $$ A(x) \cdot \left (\bm \cos 2 \pi y \\
\sin 2 \pi y \em \right )=u(x,y) \left (\bm \cos 2 \pi (y+\psi(x,y))
\\ \sin 2 \pi (y+\psi(x,y)) \em \right ). $$ The function $\psi$ is
called a {\it lift} of $A$.  Let $\mu$ be any probability measure on
$\T \times \T$ which is invariant by the continuous map $T:(x,y)
\mapsto (x+\alpha,y+\psi(x,y))$, projecting over Lebesgue measure on
the first coordinate (for instance, take $\mu$ as any accumulation
point of $\frac {1} {n} \sum_{k=0}^{n-1} T_*^k \nu$ where $\nu$ is
Lebesgue measure on $\T \times \T$). Then the number $$
\rot_f(\alpha,A)=\int \psi d\mu \mod \Z $$ does not depend on the
choices of $\psi$ and $\mu$, and is called the {\it fibered rotation
number} of $(\alpha,A)$, see \cite {JM82} and \cite {H}.

Let $$R_{\phi}=\left (\bm \cos 2\pi \phi&-\sin 2\pi \phi\\
\sin 2 \pi \phi &\cos 2\pi  \phi \em \right ),$$ then any $A\in
C^0(\T, SL(2,\R)$ is homotopic to $\theta \mapsto R_{n\theta}$ for
some $n\in\Z$, we call $n$ the degree of $A$, and denote $\deg A =n
$. The fibered rotation number is invariant under conjugation in the
following sense:  For cocycles $(\alpha,A_1)$ and $(\alpha,A_2)$, if there exists    $B \in
C^0(\T,$ $	PSL(2,\R))$, such that 
$B(\theta+\alpha)^{-1}A_1(\theta)B(\theta)=A_2(\theta),$ then we say $(\alpha,A_1)$ is conjugated to $(\alpha,A_1)$. If $B$ has degree $n$, then we have
\begin{equation}\label{rot-conj}
\rot_f(\alpha, A_1)= \rot_f(\alpha, A_2)+\frac{1}{2}  n \alpha.
\end{equation}
If furthermore $B \in C^0(\T,$ $SL(2,\R))$ with  $\deg B=n \in\Z$, then we
have
\begin{equation}\label{rot-conj'}
\rot_f(\alpha, A_1)= \rot_f(\alpha, A_2)+  n \alpha.
\end{equation}

The cocycle $(\alpha,A)$ is $C^\omega$ reducible, if it can be $C^\omega$ conjugated to a constant cocycle.
The cocycle $(\alpha,A)$ is called $C^\omega$  rotations reducible, if there exist $B \in
C^\omega(\T,SL(2,\R))$ such that
$B(\theta+\alpha)^{-1}A(\theta)B(\theta)\in SO(2,\R).$
The crucial reducibility results for us is the following:

\begin{Theorem}\cite{AFK,HoY}\label{hy1} Let
$(\alpha, A) \in \R\backslash \Q \times C^\omega_{h}(\T,SL(2,\R))$
with  $h>\tilde{h}>0,$ $R\in SL(2,\R)$, for every $\tau>1,$
$\gamma>0,$ if $\rot_f(\alpha, A)\in DC_\alpha(\tau,\gamma),$ where
$$DC_\alpha(\tau,\gamma)=\{ \phi\in \R| \|2\phi-m\alpha\|_{\R/\Z}\geq \frac{\gamma}{(|m|+1)^\tau}, m\in\Z\}$$
then there exist $T=T(\tau),$ $\kappa=\kappa(\tau)$, such that if  $$\|
 A(\theta)-R\|_h<T(\tau
)\gamma^\kappa(h-\tilde{h})^\kappa,$$ then there exist $B \in C^\omega_{\tilde{h}}(\T,SL(2,\R))$,
$\varphi \in C^\omega_{\tilde{h}}(\T,\R)$,
 such that
$$B(\theta+\alpha) A(\theta) B(\theta)^{-1}= R_{\varphi(\theta)},$$
with estimates
 $\|B-\id\|_{\tilde{h}}\leq \|
 A(\theta)-R\|_h^{\frac{1}{2}}$, 
$\|\varphi(\theta)-\hat\varphi(0)\|_{\widetilde{h}}\leq 2\|
 A(\theta)-R\|_h.$
\end{Theorem}

\subsection{Almost Mathieu cocycle}

Note that a sequence $(u_n)_{n \in \Z}$ is a formal solution of the
eigenvalue equation $H_{\lambda,\alpha,\theta} u=Eu$ if and only if
it satisfied $\begin{pmatrix}
u_{n+1}\\u_n\end{pmatrix}=S_{E}^{\lambda}(\theta+n\alpha) \cdot
\begin{pmatrix} u_n\\u_{n-1} \end{pmatrix},$ where
\begin{eqnarray*}
S_{E}^{\lambda}(\theta)=\left( \begin{array}{ccc}
 E-2\lambda\cos2\pi(\theta) &  -1\cr
  1 & 0\end{array} \right)\in SL(2,\mathbb{R}).
\end{eqnarray*}
We call  $(\alpha,S_{E}^{\lambda} )$ an almost Mathieu cocycle.

 Denote the spectrum of  $H_{\lambda,\alpha,\theta}$ by
$\Sigma_{\lambda,\alpha}$,  which  is independent of $\theta$ when
$\alpha\in \R\backslash \Q$. If $E \in \Sigma_{\lambda,\alpha}$,
then the Lyapunov exponent of almost Mathieu cocycle can be computed
directly.

\begin{Theorem}\cite{BJ}\label{bj-formula} If $\alpha \in \R
\setminus \Q$, $E \in \Sigma_{\lambda,\alpha}$, then  we have
$$L(\alpha,S_E^{\lambda})=\max
\{0,\ln |\lambda|\}.$$
\end{Theorem}

\subsection{Global theory of one frequency quasi-periodic  $SL(2,\R)$ cocycle}

We make a short review of  Avila's global theory of one frequency quasi-periodic  $SL(2,\R)$ cocycle
\cite{Aglobal}.  Suppose that  $A\in$ $C^\omega(\R/\Z,\SL(2,\R))$ admits a
holomorphic extension to $|\Im \theta|<\delta$, then for
$|\epsilon|<\delta$ we can define $A_\epsilon \in
C^\omega(\R/\Z,\SL(2,\C))$ by $A_\epsilon(\theta)=A(\theta+i \epsilon)$.
 The cocycles which are not uniformly hyperbolic are classified 
 into three regimes: subcritical, critical, and supercritical. In
 particular, $(\alpha, A)$ is said to be
 subcritical, if there exists $\delta>0,$ such that
 $L(\alpha,A_{\varepsilon})=0$ for $|\varepsilon|<\delta.$  
 
 The  heart of Avila's global theory is his \textquotedblleft Almost Reducibility Conjecture\textquotedblright (ARC), which says that subcriticality 
implies almost reducibility.
Recall the cocycle $(\alpha,A)$ is called almost reducible, if
 there exists $h_*>0$, and a sequence $B_n \in C^\omega_{h_*}(\T,PSL(2,\R))$ such that $
 B_n(\theta+\alpha)^{-1}A(\theta)B_n(\theta)$ converges to constant
 uniformly in $|\Im \theta|<h_*.$
For our purpose, we need this \textit{strong} version of almost
reducibility, and $h_*$ should be  chosen to be $\delta- \epsilon$ with $\epsilon$ arbitrary small. 


The full solution of ARC was recently given by Avila in \cite{Aac,A2}.  In the case $\beta(\alpha)>0$, it is the following:

\begin{Theorem}\cite{Aac}\label{arc}
Let $\alpha\in\R\backslash \Q$ with $\beta(\alpha)>0$, $h>0$,  $A
\in C^\omega_h(\T,\R)$. If $(\alpha, A)$ is subcritical, then for
any $0<h_*<h$ there exists $C>0$ such that if $\delta>0$ is small
enough, then there exist $B \in C^\omega_{h_*}(\T,PSL(2,\R))$ and
$R_*\in SO(2,\R)$ such that $\|B\|_{h_*}\leq e^{C\delta q}$ and
$$ \| B(\theta+\alpha)^{-1}A(\theta)B(\theta)- R_*\|_{h_*}\leq e^{-\delta q}.$$
\end{Theorem}

\subsection{Aubry duality}

Suppose that the quasi-periodic Schr\"odinger operator
\begin{equation}
(H_{V,\alpha,\theta} x)_n= x_{n+1}+x_{n-1} +V( n\alpha + \theta)
x_n=Ex_n,
\end{equation}
has an
analytic quasi-periodic Bloch wave $x_n = e^{2\pi i n\varphi}  \overline{\psi}\left(n\alpha + \phi
\right)
$ for some $
\overline{\psi}\in C^\omega(\T, \C)$ and $\varphi \in [0,1)$.  
 It is easy to see the Fourier coefficients of $\overline{\psi}(\theta)$ satisfy the following  Long-range operator:
\begin{equation}
(\widehat{L}_{V,\alpha, \varphi}u)_n=\sum _{k\in\Z}
V_ku_{n-k}+2cos2\pi (\varphi+n\alpha)u_n=Eu_n,
\end{equation}
Almost Mathieu operator is the only operator
which is invariant under Aubry duality, and the dual of
$H_{\lambda,\alpha,\theta}$ is  $H_{\lambda^{-1},\alpha,\varphi}$.

 Rigorous spectral Aubry duality was founded by
 Gordon-Jitomirskaya-Last-Simon in \cite{GJLS}, where they proved that
if $H_{\lambda,\alpha,\theta}$ has pure point spectrum for a.e.
$\theta\in\R$, then $H_{\lambda^{-1},\alpha,\varphi}$ has purely
absolutely continuous spectrum for a.e. $\varphi\in\R$. Readers can find more discussions about dynamical Aubry duality in section 4.

\section{Singular continuous spectrum}

In this section, we  prove Theorem \ref{main} (2). We re-state it as in following

\begin{Theorem}\label{singular}
Let $\alpha\in \R\backslash \Q$ with $0<\beta(\alpha)\le\infty$. If
$1\leq \lambda<e^{\beta}$, then $H_{\lambda,\theta,\alpha}$ has
purely singular continuous spectrum for any $\theta\in\T$.
\end{Theorem}

\begin{Remark}\label{gordon}
We stress again by classical Gordon's argument \cite{G}, one can only 
obtain result in rigime $1\leq \lambda<e^{\frac{\beta}{2}}$.  The
reason why one can only obtain  $e^{\frac{\beta}{2}}$ is that, in the
classical Gordon's lemma,  one has to approximate the
solution by periodic ones  along double periods. 
\end{Remark}

\begin{pf}
If $1<\lambda<e^{\beta}$, $E \in \Sigma_{\lambda,\alpha}$, then by
Theorem \ref{bj-formula}, one always has $L(E,\alpha)=\ln\lambda>0$.
By Kotani's theory \cite{Ko84}, the operator
$H_{\lambda,\theta,\alpha}$ doesn't support any absolutely
continuous spectrum, thus one only needs to exclude the point
spectrum. In the case  $\lambda=1$, since Lebesgue measure of
$\Sigma_{1,\alpha}$ is zero for any  $\alpha\in \R\backslash \Q$
\cite{AK06,L94}, then $H_{1,\theta,\alpha}$ also doesn't support any
absolutely continuous spectrum, thus to prove Theorem
\ref{singular}, it is also enough to exclude the point spectrum.

As in classical Gordon's lemma, we approximate the quasi-periodic
cocycles by periodic ones. Denote
$A(\theta)=S_{E}^{\lambda}(\theta)$ and
\begin{eqnarray}A_m(\theta)&=&A(\theta+(m-1)\alpha)\cdots A(\theta+\alpha)A(\theta),\\
                           \nonumber &=&A^m(\theta)\cdots A^2(\theta)A^1(\theta) \end{eqnarray}
 \begin{eqnarray}\widetilde{A}_m(\theta)&=&A(\theta+(m-1)\frac{p_n}{q_n})\cdots
A(\theta+\frac{p_n}{q_n})A(\theta),\\
\nonumber  &=&\widetilde{A}^m(\theta)\cdots
\widetilde{A}^2(\theta)\widetilde{A}^1(\theta) , \end{eqnarray} for $m\geq1.$ 
We also denote $A_{-m}(\theta)=A_m(\theta-m\alpha)^{-1},$
$\widetilde{A}_{-m}(\theta)=\widetilde{A}_m(\theta-m\frac{p_n}{q_n})^{-1}.$
Our proof is based on the following

\begin{Proposition}\label{appro}
Let $\alpha\in \R\backslash \Q$. If
$\lambda \geq 1$ and $E\in \Sigma_{\lambda,\alpha}$, then
for any $\epsilon>0$, there exists $N=N(E,\lambda, \epsilon)>0$ such
that if $q_n>N$, then we have
\begin{eqnarray}
\label{appro-2}\sup_{\theta\in\T}\|\widetilde{A}_{ \pm
q_n}(\theta)-A_{\pm q_n}(\theta)\| & \leq& \frac{1}{q_{n+1}}
e^{(\ln\lambda+ \epsilon)q_n},\\
\label{appro-7}\sup_{\theta\in\T}\| A_{q_n}(\theta+q_n\alpha)-
A_{q_n}(\theta)\|& \leq& \frac{1}{q_{n+1}}
e^{(\ln\lambda+ \epsilon)q_n}.
\end{eqnarray}

\end{Proposition}

\begin{pf}
Furman's result \cite{F} gives
\begin{eqnarray}\label{appro-3}\lim_{m \rightarrow
\pm \infty} \sup_{\theta \in \T} \frac{1}{|m|}\log
\|A_m(\theta)\|\le L(\alpha,S_E^{\lambda}).\end{eqnarray}
    Then by Theorem
    \ref{bj-formula}, we know for any $\epsilon>0$, there exists
$K=K(E,\lambda,\epsilon)>0$, such that for any $|m| \geq K$, we have
\begin{eqnarray}\label{appro-6}
\sup_{\theta\in\T}\|A_{m}(\theta)\| & \leq&
e^{|m|(\ln\lambda+\epsilon/2)}.
\end{eqnarray}

In the following, we only consider $m$ is positive,  the proof is similar for negative $m$. In order to prove $(\ref{appro-2}),$
we need the following:
\begin{Lemma}\label{appro-lemma}
Let $\alpha\in \R\backslash \Q$. If $
\lambda \geq 1$ and $E\in \Sigma_{\lambda,\alpha}$, then for any
$\epsilon>0$, there exists $N_-=N_-(E,\lambda, \epsilon)>2K$, such
that
\begin{eqnarray} \label {3.7}
\sup_{\theta\in\T}\|\widetilde{A}_{m}(\theta)\| & \leq&
e^{m(\ln\lambda+2 \epsilon/3)}
\end{eqnarray}for any $q_n \geq N_-$,  $m \geq K$.
\end{Lemma}

\begin{pf}

Clearly, for fixed $m \in \Z$ and $\delta>0$, if $q_n$ is
sufficiently large we have
$$\sup_{\theta\in\T}\big|\frac{1}{m} \ln  \|\widetilde{A}_m(\theta)\|- \frac{1}{m} \ln  \|A_m(\theta)\| \big|<\delta.$$
Thus, there exists $N_-=N_-(E,\lambda,\epsilon)>0$ such that if $q_n
\geq N_-$ then $(\ref {3.7})$ holds for $K \leq m \leq 2K-1$. Since
any $m \geq K$ can be written as a sum of integers $m_i$ satisfying
$K \leq m_i \leq 2K-1$, this implies that $(\ref {3.7})$ holds for
all $m \geq K$.
\end{pf}

Once we have Lemma \ref{appro-lemma}, $(\ref{appro-2})$ can be
proved directly by telescoping arguments. In fact, if $q_n \geq N_-$
we can write
\begin{eqnarray*}
A_{q_n} - \widetilde{A}_{q_n} &=&  \sum_{i=1}^{q_n}A^{q_n}\cdots
A^{i+1}\Big(A^{i}-\widetilde{A}^{i}\Big)\widetilde{A}^{i-1}\cdots
\widetilde{A}_1\\ &=& \Big(\sum_{i=1}^{K}+ \sum_{i=K+1}^{q_n-K}
+\sum_{i=q_n-K+1}^{q_n}\Big)A^{q_n}\cdots
A^{i+1}\Big(A^{i}-\widetilde{A}^{i}\Big)\widetilde{A}^{i-1}\cdots
\widetilde{A}_1\\
&=& (I)+(II)+(III),
\end{eqnarray*}
since for $i\leq q_n$, we have $\|A^{i}-\widetilde{A}^{i}\|\leq
\frac{4 \pi \lambda (i-1)}{q_{n}q_{n+1}}\leq \frac{4 \pi
\lambda}{q_{n+1}},$ then by $(\ref{appro-6})$ and Lemma
\ref{appro-lemma}, we can estimate
\begin{eqnarray*}
(I)&\leq& \frac{4 \pi \lambda}{q_{n+1}} \sum_{i=1}^{K}
(4\lambda+3)^{i-1}e^{(q_n-i)(\ln\lambda+2\epsilon/3)}, \\
(II)&\leq& \frac{4 \pi \lambda}{q_{n+1}}
\sum_{i=K+1}^{q_n-K}e^{(q_n-1)(\ln\lambda+2\epsilon/3)},\\
(III)&\leq & \frac{4 \pi \lambda}{q_{n+1}} \sum_{i=q_n-K+1}^{q_n}
(4\lambda+3)^{q_n-i}e^{(i-1)(\ln\lambda+2\epsilon/3)}.
\end{eqnarray*}
If $q_n$ is sufficiently large, then $(\ref{appro-2})$ follows
directly. Using the similar argument as above, we can prove $(\ref{appro-7})$.

\end{pf}

Now we finish the proof of Theorem \ref{singular} by contradiction.
For any fixed $\theta$, we suppose that  $E$ is an eigenvalue of
$H_{\lambda,\alpha,\theta}$, then there exists  $\overline{v}= \begin{pmatrix} v_0\\v_{-1}
\end{pmatrix} $ with  $\|\overline{v}\|=1,$ and  for  any
$\varepsilon>0$, there exists
$\overline{N}=\overline{N}(E,\lambda,\varepsilon)$, such that if
$|m|> \overline{N}(E,\lambda,\varepsilon)$, then
$\|A_m(\theta)\overline{v}\|\leq \varepsilon.$
In particular,  for any $0<2\epsilon< \ln\lambda-\beta$,  we can
select $q_n>\max\{ N(E,\lambda,\epsilon),$ $
\overline{N}(E,\lambda,\varepsilon)\}$, and
$q_{n+1}>e^{(\beta-\epsilon)q_n}$, such that
\begin{equation}\label{initial} \|A_{q_n}(\theta)\overline{v}\|\leq
\varepsilon,
\qquad \|A_{-{q_n}}(\theta)\overline{v}\|\leq \varepsilon,
\end{equation} where
$N(E,\lambda,\epsilon)$ is defined in Proposition \ref{appro}.

What's important is the following observation:
\begin{Lemma}\label{trace}The following estimate holds:
\begin{eqnarray*}
\|A_{q_n}(\theta+q_n\alpha)+A_{-{q_n}}(\theta+q_n\alpha)\| \leq 2\varepsilon+ 10 e^{-(\beta-\ln\lambda-2
\epsilon)q_n}.
\end{eqnarray*}
\end{Lemma}

\begin{pf}
By $(\ref{appro-2})$, it is sufficient for us to prove
\begin{eqnarray}\label{midesti}
\|\widetilde{A}_{q_n}(\theta+q_n\alpha)+\widetilde{A}_{-{q_n}}(\theta+q_n\alpha)\| \leq 2\varepsilon+ 8 e^{-(\beta-\ln\lambda-2
\epsilon)q_n}.
\end{eqnarray}
By  Hamilton-Clay Theorem, for any $M\in SL(2,\R)$, one has
  \begin{eqnarray}\label{hamicaly}
  M+M^{-1}=\tr M\cdot Id,
  \end{eqnarray}
  for every $\theta' \in \T$. Take
$M=\widetilde{A}_{{q_n}}(\theta')$, then
\begin{eqnarray}\label{hami-caly}
\widetilde{A}_{{q_n}}(\theta')+\widetilde{A}_{-{q_n}}(\theta')=\tr
\widetilde{A}_{q_n}(\theta').\end{eqnarray}

By assumptions $(\ref{initial})$ and $(\ref{appro-2})$, we have
\begin{eqnarray*}&&\| \tr
\widetilde{A}_{q_n}(\theta)\|\\
\nonumber &\leq& \| A_{q_n}(\theta)\overline{v}+ A_{-{q_n}}(\theta)
\overline{v}\|+\|\widetilde{A}_{q_n}(\theta)-A_{q_n}(\theta)\|+\|\widetilde{A}_{-{q_n}}(\theta)-A_{-{q_n}}(\theta)\|\\
\nonumber &\leq& 2\varepsilon+2 e^{-(\beta-\ln\lambda-2
\epsilon)q_n}.
\end{eqnarray*}
As a result of $(\ref{appro-2})$ and $(\ref{appro-7})$, we have
\begin{eqnarray*}
&&\| \tr \widetilde{A}_{q_n}(\theta+q_n\alpha)\|\\ &\leq& \| \tr
\widetilde{A}_{q_n}(\theta+q_n\alpha)- \tr
A_{q_n}(\theta+q_n\alpha)\|+\|\tr A_{q_n}(\theta+q_n\alpha)- \tr
A_{q_n}(\theta)\|\\
&& +\|\tr A_{q_n}(\theta)-\tr \widetilde{A}_{q_n}(\theta)\|+\|\tr
\widetilde{A}_{q_n}(\theta)\|\\
&\leq & 2\varepsilon+8 e^{-(\beta-\ln\lambda-2 \epsilon)q_n},
\end{eqnarray*}
then $(\ref{midesti})$ follows from $(\ref{hami-caly})$.
\end{pf}

However  by Lemma \ref{trace}, we have 
\begin{eqnarray*}
&&\|A_{2q_n}(\theta)\overline{v}\|=
\|A_{q_n}(\theta+q_n\alpha) A_{q_n}(\theta)\overline{v}  \|\\
&\geq& \|A_{-{q_n}}(\theta+q_n\alpha) A_{q_n}(\theta)\overline{v} \|- \|\widetilde{A}_{q_n}(\theta+q_n\alpha)+\widetilde{A}_{-{q_n}}(\theta+q_n\alpha)\| \|A_{q_n}(\theta)\overline{v}\| \\
&\geq&1- 2\varepsilon^2-10 \varepsilon e^{-(\beta-\ln\lambda-2
\epsilon)q_n}> \frac{1}{2},\end{eqnarray*} which contradicts with
the assumption that $E$ is an eigenvalue.

\end{pf}

%
%

\section{Anderson localization}\label{ander}

In this section, we  prove Theorem \ref{main} (3). We re-state it as the following
\begin{Theorem}\label{anderson-transition}
Let $\alpha\in\R\backslash \Q$ be such that $0<\beta(\alpha)<\infty.$
If $\lambda>e^{\beta},$  then the almost Mathieu operator
$H_{\lambda,\alpha,\phi}$ has Anderson Localization for a.e. $\phi$. 
\end{Theorem}

Traditional method for Anderson Localization  is to prove the exponentially decay of Green function \cite{AJ05,J94,J95,J99}. Due to the limitation of the method,  Anderson Localization can be proved only for Liouvillean frequency with $\lambda>e^{ \frac{16\beta}{9}}$ so far  \cite{AJ05}. So there is still a gap between $e^{\beta}$ and $e^{ \frac{16\beta}{9}}$.

In this  paper, we develop a new approach  depending on the reducibility and  Aubry duality.
 We will show that Theorem \ref{anderson-transition} can be obtained by 
dynamical Aubry duality and the following full measure reducibility result:
 \begin{Theorem}\label{full}
Let $\alpha \in \R \setminus \Q$ with $\beta(\alpha)>0$,  if
$\lambda>e^\beta$, $\rot_f(\alpha, S_E^{\lambda^{-1}})$ is
Diophantine w.r.t. $\alpha$, then $(\alpha, S_E^{\lambda^{-1}})$ is
reducible. 
\end{Theorem}

The dynamical Aubry duality was established by Puig \cite{Pui06},
who proved that Anderson
 localization of the Long range operator $\widehat{L}_{V,\alpha, \varphi}$ for almost every $\varphi\in\T$ implies reducibility of $(\alpha,S_E^V)$ for almost every
energies.    Conversely,  to deal with localization problem by reducibility was first realized by 
  You-Zhou in \cite{YZ}. However, in \cite{YZ} they can only prove the eigenvalues of $\widehat{L}_{V,\alpha, \varphi}$ with exponentially decaying eigenfunctions are dense in the spectrum. The main issue remained is to prove  those eigenfunctions form a complete basis. The key point in this paper is that, we find  that the quantitative estimates in the proof of Theorem \ref{full} actually provides 
an asymptotical
distribution of the eigenvalues and eigenfunctions, which ultimately implies  pure point spectrum for almost every phases. 
Compared with tradition localization argument, the price we have to pay is that we lose precise arithmetic  control on  the localization phases.
However, by this approach, one can indeed  establish a kind of  equivalence
between quantitative full measure reducibility of Schr\"odinger
operator (or Schr\"odinger cocycle) and Anderson
 localization of its dual Long-range operator.\\

\noindent
\textbf{Proof Theorem \ref{anderson-transition}:}   We
need the following definition:

\begin{Definition}
For any fixed $N\in\N,C>0,\varepsilon>0$, a normalized eigenfunction
$u(n)$ is said to be  $(N,C,\varepsilon)$-good, if $|u(n)|\leq
e^{-C\varepsilon|n|}$, for $|n|\geq (1-\varepsilon)N$.
\end{Definition}

We label the $(N,C,\varepsilon)$-good eigenfunctions of $H_{\lambda,
\alpha,\phi}$ by $u_j^\phi(n)$,  denote the corresponding
eigenvalue by $E_j^\phi$, also we  denote
$$ \mathcal {E}_{N,C,\varepsilon}^{\phi}=\{E_j^\phi| u_j^\phi(n) \text{ is a
$(N,C,\varepsilon)$-good normalized eigenfunction}\}$$
and denote $\mathcal {E}(\phi)= \bigcup_{N>0}  \mathcal {E}_{N,C,\varepsilon}^{\phi}.$
 Let $\mu_{\delta_0,\phi}^{pp}$ be the spectral measure supported on
$\mathcal {E}(\phi)$ with respect to $\delta_0$.\\

The following spectral analysis is completely new and  will be crucial for our proof.
\begin{Proposition}\label{distribution}
Suppose that there exists $C>0$, such that for any $\delta>0,$ there
exists $\varepsilon>0$, and for a.e. $\phi$,
\begin{equation}\label{good} \#\{\text{linearly independent
$(N,C,\varepsilon)$-good eigenfunctions}\}\geq
(1-\delta)2N,\end{equation} for $N$ large enough, then for a.e.
$\phi$, we have  $\mu_\phi=\mu_{\delta_0,\phi}=
\mu_{\delta_0,\phi}^{pp}$.
\end{Proposition}

\begin{pf}

Fix $\phi\in\T^1$ such that (\ref{good}) is satisfied. Denote
$$ K_{N,C,\varepsilon}^{\phi}=\{ j\in\N| u_j^\phi(n) \text{ is a
$(N,C,\varepsilon)$-good eigenfunction}\}.$$ Notice that for any fixed
${N,C,\varepsilon}$, $\# K_{N,C,\varepsilon}^{\phi}$ is finite, and also \begin{equation}\label{eigen}
\sum_{|n|\leq (1-\varepsilon)N}|u_j^\phi(n)|^2 >1- e^{-C\varepsilon
N},
\end{equation}
for  $(N,C,\varepsilon)$-good
eigenfunction $u_j^\phi(n)$.
 
 Let $\widetilde{\mu}^{pp}_{\delta_n,\phi}=
\widetilde{\mu}^{pp}_{\delta_n,\phi}(N,C,\varepsilon)$ be the
truncated spectral measure supported on $\mathcal
{E}_{N,C,\varepsilon}^{\phi}$. Then by spectral theorem and
$(\ref{eigen})$, we have
\begin{eqnarray*}
\frac{1}{2N}\sum_{|n|\leq N}|\mu^{pp}_{\delta_n,\phi}|&>&
\frac{1}{2N}\sum_{|n|\leq N}|\widetilde{\mu}^{pp}_{\delta_n,\phi}|\\
&=&\frac{1}{2N}\sum_{|n|\leq N}\langle P_{\mathcal
{E}_{N,C,\varepsilon}^{\phi}}\delta_n,
\delta_n\rangle\\
&=&\frac{1}{2N}\sum_{|n|\leq N}\sum_{j\in
K_{N,C,\varepsilon}^{\phi}}\langle P_{E_j^\phi}\delta_n,
\delta_n\rangle\\
&>&\frac{1}{2N}\sum_{|n|\leq(1-\varepsilon)N}\sum_{j\in
K_{N,C,\varepsilon}^{\phi}}|u_j^\phi(n)|^2\\
&>&\frac{1}{2N}  \# K_{N,C,\varepsilon}^{\phi}(1- e^{-C\varepsilon N})\\
&>& (1-\delta)(1-e^{-C\varepsilon N}).
\end{eqnarray*}
Since $ \mathcal {E}(\phi) =\mathcal {E}(\phi+\alpha)$, 
we can rewrite the above inequalities as $$ \frac{1}{2N}\sum_{|n|\leq N}|\mu^{pp}_{\delta_0,\phi+n\alpha}|>  (1-\delta)(1-e^{-C\varepsilon N}),$$
Let $N$ go to $\infty$, since $\delta$ is arbitrary small,  we have
$$ \int_{\T^1}
| \mu^{pp}_{\delta_0,\phi}|d\phi=1,$$ by Birkhoff's ergodic theorem. 
 Thus for $a.e.$ $\phi\in\T^1$,
$\mu_\phi=\mu_{\delta_0,\phi}= \mu_{\delta_0,\phi}^{pp}.$
\end{pf}

Let $\Theta_\gamma=\{\phi| \phi\in DC_\alpha(\tau,\gamma)
\}$. We have $ \bigcup_{\gamma>0}\Theta_\gamma=1$, which
implies that for any $\delta>0$, there exists $
\widetilde{\varepsilon}>0$, such that if
$|\gamma|<\widetilde{\varepsilon},$ then $|\Theta_\gamma |>
1-\frac{\delta}{3}.$ By Birkhoff's ergodic theorem again, we have
$$\lim_{\widetilde{N} \rightarrow \infty} \frac{1}{2\widetilde{N}} \sum_{|k|\leq \widetilde{N}} \chi_{\Theta_\gamma }(\phi+k\alpha)= \int_{\T^1}\chi_{\Theta_\gamma }(\phi)d\phi. $$
Thus for $N$ large enough (we take
$\widetilde{N}=N(1-\frac{\delta}{3})$), we have
\begin{equation}\label{number}
\#\{k| \phi+k\alpha \in \Theta_\gamma, |k|\leq
2N(1-\frac{\delta}{3}) \}\geq (1-\delta)2N.\end{equation}

For any $\phi \in \Theta_\gamma$, we choose $\bar N$ sufficiently
large such that (\ref{number}) holds for $N>\bar N$. We will prove
that $H_{\lambda, \alpha,\phi}$ has at least  $(1-\delta)2N$ different
eigenvalues $E_k^\phi$ whose eigenfunctions $u_k^\phi(n)$  are $(N,\ln
\lambda -\beta-\epsilon, \varepsilon)$-good for any $\epsilon$. To
prove this, we need the following \textit{quantitative} version of Theorem \ref{full}:

\begin{Proposition}\label{prop}
Let $\alpha \in \R \setminus \Q$ with $\beta(\alpha)>0$ and
$\lambda>e^\beta$.  Suppose that $\rot_f(\alpha,
S_{\lambda^{-1}E_k}^{\lambda^{-1}})=\phi+k\alpha \in
DC_\alpha(\tau,\gamma)$. Then  for any fixed $\gamma>0$, $\tau>0$
and small enough $\epsilon>0$,  there exist
$c_1(\lambda, \gamma,\tau,\epsilon,\alpha), c_2(\lambda,\gamma,\tau,\epsilon)$  and $B_k \in
C^\omega_{\ln\lambda-\beta-\epsilon}(\T,SL(2,\R))$, such that
\begin{equation}\label{prop-1} B_k(\theta+\alpha)^{-1}
S_{\lambda^{-1}E_k}^{\lambda^{-1}}(\theta)B_k(\theta)=R_{\phi+k^{'}\alpha},\end{equation}
with estimates:
\begin{eqnarray}
\label{esti-1} \| B_k\|_{\ln\lambda-\beta-\epsilon} &\leq&
c_1(\lambda,\gamma,\tau,\epsilon,\alpha),\\
 \label{esti-2} |k-k^{'}|&\leq& c_2(\lambda,\gamma,\tau,\epsilon).
\end{eqnarray}

\end{Proposition}

\begin{pf}

If $\lambda>e^{\beta}>1$, $\lambda^{-1} E_k \in \Sigma_{\lambda^{-1},\alpha}$, then the almost Mathieu cocycle $(\alpha,
S_{\lambda^{-1}E_k}^{\lambda^{-1}})$ is subcritical in the regime
$|\mathfrak{I}\theta|<\ln\lambda$.  To prove Proposition \ref{prop}, we  need 
 Theorem \ref{bj-formula} and the following:

\begin{Lemma}
If $\alpha \in \R \setminus \Q$, $\lambda>1$, $E \in \R$, then for
$\epsilon \geq 0$,
$$L(\alpha,(S_E^{\lambda^{-1}})_\epsilon)=\max
\{L(\alpha,S_E^{\lambda^{-1}}),(\epsilon-\ln \lambda)\}.$$
\end{Lemma}

\begin{pf} The proof can be found in Appendix A of  \cite{Aglobal}.
\end{pf}

Now by Theorem \ref{arc}, for $0<2\epsilon<\ln \lambda-\beta$,  there exists a sequence of
$\widetilde{B}_n \in
C^\omega_{\ln\lambda-\epsilon/2}(\T,PSL(2,\R))$  such that
$$\widetilde{B}_n(\theta+\alpha)^{-1}S_{\lambda^{-1}E_k}^{\lambda^{-1}}(\theta)\widetilde{B}_n(\theta)= R_{\varphi_n}+F_n(\theta),$$
with estimate
\begin{eqnarray}
\label{esti-1'}\|\widetilde{B}_n\|_{\ln\lambda-\epsilon/2}&\leq&
e^{C\delta^{'} q_n},\\
\nonumber \|F_n\|_{\ln\lambda-\epsilon/2}&\leq&
e^{-\delta^{'} q_n},
\end{eqnarray}
which implies  
\begin{equation}\label{deg1}
|\deg \widetilde{B}_n| \leq c(\lambda,\epsilon) q_n.
\end{equation}
One may consult footnote 5 of \cite{Aac} in proving this.

If $\phi+k\alpha \in DC_\alpha(\tau,\gamma),$ we have 
\begin{eqnarray*}
&&\|2(\phi+k\alpha)-m\alpha-k'\alpha\|_{\R/\Z}\\
&\geq& \frac{\gamma}{(|m+k^{'}|+1)^\tau} \geq
\frac{(1+|k^{'}|)^{-\tau}\gamma}{(|m|+1)^\tau}.
\end{eqnarray*}
By   $(\ref{rot-conj})$, this formula  implies that  $ \rot_f(\alpha, R_{\varphi_n}+F_n(\theta)) \in
DC_\alpha(\tau,(1+|\deg \widetilde{B}_n|)^{-\tau}\gamma)$.  Let  $q_s$  be the smallest denominator such that
\begin{eqnarray*}
q_{s+1}&> &e^{(\beta-o(1))q_s},\\
e^{-q_s \delta^{'}} &<& T(\tau)(\frac{\gamma}{(1+c(\lambda,\epsilon)|q_s|)^{\tau}})^\kappa(\frac{\epsilon}{2})^\kappa,
\end{eqnarray*}
where $T=T(\tau),$ $\kappa=\kappa(\tau)$ are defined in Theorem
\ref{hy1}. By Theorem \ref{hy1},  there
exist $\overline{B}_k(\theta) \in
C^\omega_{\ln\lambda-\epsilon}(\T,SL(2,\R)),$ $\eta_k(\theta) \in
C^\omega_{\ln\lambda-\epsilon}(\T,\R),$ such that 
$$\overline{B}_k(\theta+\alpha)^{-1}(  R_{\varphi_s}+F_s(\theta))\overline{B}_k(\theta)=R_{\eta_k(\theta)}.$$
with estimates $ \| \eta_k  \|_{\ln\lambda-\epsilon} \leq e^{-q_s \delta^{'}} $ and
\begin{equation} \label{esti-9}
\| \overline{B}_k-id  \|_{\ln\lambda-\epsilon} \leq e^{-q_s \delta^{'}/2} .
\end{equation} Let $\psi_k(\theta)$ satisfy
\begin{equation}\label{homo}
\psi_k(\theta+\alpha)-\psi_k(\theta)=\eta_k(\theta)-\hat{\eta}_k(0).
\end{equation}
since  $\ln \lambda>\beta$, by $(\ref{equibeta})$, we know that there exists $c=c(\alpha,\epsilon)$ such that $(\ref{homo})$ has analytic solution $\psi_k(\theta) \in
C^\omega_{\ln\lambda-\beta-\epsilon}(\T,\R)$ with estimate 
\begin{equation}\label{esti-10}
\|\psi_k\|_{\ln\lambda-\beta-\epsilon} \leq c(\alpha,\epsilon)  \| \eta_k  \|_{\ln\lambda-\epsilon} \leq  c(\alpha,\epsilon) e^{-q_s \delta^{'}} .
\end{equation}
Let $B_k(\theta)= \widetilde{B}_s(\theta)\overline{B}_k(\theta)R_{\psi_k(\theta)}$, then there exists $k^{'}\in \Z$, such that 
$$B_k(\theta+\alpha)^{-1}
S_{\lambda^{-1}E_k}^{\lambda^{-1}}(\theta)B_k(\theta)=R_{\hat{\eta}_k(0)}=  R_{\phi+k^{'}\alpha}.$$
Since $\rot_f(\alpha,
S_{\lambda^{-1}E_k}^{\lambda^{-1}})$ is irrational w.r.t $\alpha$, then  $B_k(\theta) \in
C^\omega_{\ln\lambda-\beta-\epsilon}(\T,SL(2,\R))$,  one can consult  Remark 1.5 of \cite{AK06} for this proof.  Notice that   $\deg R_{\psi_k(\theta)}=0$ and  by $(\ref{esti-9})$, we have  $\deg \overline{B}_k =0$.  Consequently by $(\ref{rot-conj'})$,
 we have
\begin{eqnarray}
\label{deg2}k^{'}=k- \deg \widetilde{B}_s.\end{eqnarray} 
$(\ref{esti-2})$ then follows from $(\ref{deg1})$ and $(\ref{deg2})$, and
$(\ref{esti-1})$ follows from $(\ref{esti-1'})$, $(\ref{esti-9})$ and $(\ref{esti-10})$.
\end{pf}

Rewrite $(\ref{prop-1})$ as
\begin{equation}\label{a-1} B_k(\theta+\alpha)^{-1}
S_{\lambda^{-1}E_k}^{\lambda^{-1}}(\theta)B_k(\theta)= \left(
\begin{array}{ccc}
e^{2\pi i(\phi+k^{'}\alpha)}& 0\cr
  0 & e^{-2\pi i(\phi+k^{'}\alpha)}\end{array} \right),\end{equation}
and write $B_k(\theta)=\left(
\begin{array}{ccc}
z_{11}(\theta) &  z_{12}(\theta) \cr z_{21}(\theta) &z_{22}(\theta)
\end{array} \right),$ then
we have
\begin{eqnarray}\label{block-red}&& (\lambda^{-1}E_k-2\lambda^{-1}
\cos(\theta))z_{11}(\theta)\\ \nonumber&=&
z_{11}(\theta-\alpha)e^{-2\pi
i(\phi+k^{'}\alpha)}+z_{11}(\theta+\alpha)e^{2\pi
i(\phi+k^{'}\alpha)}.\end{eqnarray} Taking the Fourier
transformation for $(\ref{block-red})$, we have
\begin{eqnarray*}
\widehat{z}_{11}(n+1)+\widehat{z}_{11}(n-1)+2\lambda\cos2\pi
(\phi+k^{'}\alpha+n\alpha)\widehat{z}_{11}(n)=
E_k\widehat{z}_{11}(n),
\end{eqnarray*}
then $\widehat{z}_{11}(n)$ is a eigenfunction, since $z_{11}\in
C^\omega_{\ln\lambda-\beta-\epsilon}(\T,\C)$. To normalize
$\widehat{z}_{11}(n)$, we need the following observation:

\begin{Lemma}\label{z1-estimate}
We have the following:
 $$\|\widehat{z}_{11}\|_{l^2}\geq (2\|B\|_{C^0})^{-1}.$$
\end{Lemma}

\begin{pf}
Write $$u=\left(\begin{array}{c} z_{11}(\theta) \\
z_{21}(\theta)
\end{array}\right), \qquad v=\left(\begin{array}{c} z_{12}(\theta) \\
z_{22}(\theta)
\end{array}\right),$$ then
$\|u\|_{L^2}\|v\|_{L^2}>1$ since $\det B_k(\theta)=1.$ This implies that 
$$ \|z_{11}\|_{L^2}+ \|z_{21}\|_{L^2}=  \|u\|_{L^2}> \|v\|_{L^2}^{-1}>(\|B\|_{C^0})^{-1}.$$
 By $(\ref{a-1})$, we have
$z_{21}(\theta+\alpha)=e^{-2\pi i(\phi+k^{'}\alpha)}z_{11}(\theta),$
therefore, we have
$$\|\widehat{z}_{11}\|_{l^2}=\|z_{11}\|_{L^2}  \geq (2\|B\|_{C^0})^{-1}.$$
\end{pf}

Normalizing $\widehat{z}_{11}(n)$ by $u_k^{\phi}(n)=\frac{\widehat{z}_{11}(n+k^{'})}{
\|\widehat{z}_{11}\|_{l^2}}$.
Now we prove it is in fact $(N,\ln \lambda -\beta-\epsilon,
\varepsilon)$-good. Let
$$  2\varepsilon< \frac{\delta}{3}- \frac{c_3(\lambda,\gamma,\tau,\epsilon,\alpha)}{N}- \frac{c_2(\lambda,\gamma,\tau,\epsilon)}{N},
$$ where $c_3(\lambda, \gamma,\tau,\epsilon,\alpha)=\frac{\ln 2c_1(\lambda,\gamma,\tau,\epsilon,\alpha)}{ \ln\lambda-\beta-\epsilon }.$ Since
$u_k^{\phi}(n)=u_k^{\phi+k^{'}\alpha}(n-k^{'})$,  then by Proposition \ref{prop} and Lemma \ref{z1-estimate},
 we
have
\begin{eqnarray*}
|u_k^{\phi}(n) | &=&|u_k^{\phi+k^{'}\alpha}(n-k^{'})|\\
&\leq&  \| B_k\|_{\ln\lambda-\beta-\epsilon}^2
e^{-|n-k^{'}| (\ln \lambda -\beta-\epsilon)}\\
&\leq & e^{ ( c_3(\lambda,\gamma,\tau,\epsilon,\alpha)+
|k|+c_2(\lambda,\gamma,\tau,\epsilon))(\ln \lambda
-\beta-\epsilon)} e^{-|n| (\ln \lambda -\beta-\epsilon)}\\
&\leq&
 e^{ (N(1-\frac{\delta}{3})+c_2(\lambda,\gamma,\tau,\epsilon)+c_3(\lambda,\gamma,\tau,\epsilon,\alpha))(\ln
\lambda -\beta-\epsilon)}e^{-|n|(\ln \lambda -\beta-\epsilon) }\\
&\leq& e^{-|n|(\ln \lambda -\beta-\epsilon)\varepsilon},
\end{eqnarray*}
for  $|n|\geq N(1-\varepsilon)$, which means $ (u_k^{\phi}(n))$ is
$(N,\ln \lambda -\beta-\epsilon, \varepsilon)$-good.

By Proposition \ref{distribution} and the above estimate, we know for  $a.e.$
$\phi\in\T^1$,
$H_{\lambda,\alpha,\phi}$ has Anderson Localization.

\section*{Acknowledgements}
A.A was partially supported by the ERC Starting Grant\textquotedblleft Quasiperiodic\textquotedblright and
by the Balzan project of Jacob Palis. J. Y was partially supported by NNSF of China (11471155) and
973 projects of China (2014CB340701).   Q. Z was partially supported by
Fondation Sciences Math\'{e}matiques de Paris (FSMP) and  and ERC
Starting Grant \textquotedblleft Quasiperiodic\textquotedblright.

\end{document}